\documentclass[12pt]{amsart}
\usepackage{geometry}                
\usepackage{graphicx}
\usepackage{amsrefs}
\usepackage{amssymb}
\usepackage[matrix,arrow,ps]{xy}
\usepackage{epstopdf}
\newtheorem{teorema}{Theorem}[section]
\newtheorem{lemma}[teorema]{Lemma}
\newtheorem{propos}[teorema]{Proposition}
\newtheorem{corol}[teorema]{Corollary}
\newtheorem{ex}{Example}[section]
\newtheorem{rem}{Remark}[section]
\newtheorem{defin}[teorema]{Definition}
\def\bt{\begin{teorema}}
\def\et{\end{teorema}}
\def\bp{\begin{propos}}
\def\ep{\end{propos}}
\def\bl{\begin{lemma}}
\def\el{\end{lemma}}
\def\bc{\begin{corol}}
\def\ec{\end{corol}}
\def\br{\begin{rem}\rm}
\def\er{\end{rem}}
\def\bex{\begin{ex}\rm}
\def\eex{\end{ex}}
\def\bd{\begin{defin}}
\def\ed{\end{defin}}
\def\demo{\par\noindent{\bf Proof.\ }}
\def\enddemo{\ $\Box$\par\vskip.6truecm}
\def\R{{\mathbb R}}   \def\a {\alpha} \def\b {\beta}
\def\N{{\mathbb N}}      \def\e{\varepsilon}
\def\C{{\mathbb C}}      \def\l{\lambda}
                 \def\o{\omega}\def\r{\varrho}
                   
\def\P{{\mathbb P}}

 \def\oli{\overline}

\def\O{\Omega}

\title {Non-embeddability of certain classes of Levi flat manifolds}

\author{Giuseppe Della Sala}
\address{Department of Mathematics\\ \newline University of Vienna\\ Vienna\\ 1090\\ Austria}
\email{giuseppe.dellasala@univie.ac.at, beppe.dellasala@gmail.com}

\begin{document}

\begin{abstract}
On the basis of a result of Barrett \cite{Ba}, we show that members of certain classes of abstract Levi flat manifolds with boundary, whose Levi foliation contains a compact leaf with contracting, flat holonomy, admit no $CR$ embedding as a hypersurface of a complex manifold.
\end{abstract}

\maketitle

In \cite{Ba}, Barrett showed that there is no Levi flat submanifold $S\cong S^3$, smoothly embedded in a complex $2$-manifold $M$, such that its foliation is diffeomorphic to Reeb's one. A key ingredient in the proof is a result by Ueda \cite{Ue}, which allows to find an equation for a compact complex curve $C\subset M$ (in a neighborhood of $C$), provided that its normal bundle satisfies certain triviality conditions.

We show that Barrett's method can be adapted to prove that other classes of Levi flat manifolds, of dimension greater than $3$, are non-embeddable as smooth hypersurfaces of a complex manifold. This is due to the fact that the relevant part of Ueda's argument is valid also in dimension greater than $1$.

In our situation, we assume the existence of a compact leaf whose holonomy is contracting and flat, as in the case of Reeb's foliation. Moreover, we ask for the holonomy covering of the compact leaf to be (partially) \lq\lq extendable\rq\rq\ at infinity, a technical condition (based on the notion of partial compactification employed in \cite{MV2}) which is very often verified in practice - as in Reeb's case, and in the case of the examples discussed in section \ref{second}.

The proof of theorem \ref{main} is pretty much contained in \cite{Ba}, \cite{MV2} and \cite{Ue}; the purpose of this note is essentially to explain it in detail, and a part of our argument is in fact pointing out why Theorem 3 in \cite{Ue} applies to our situation. Once a defining function for the compact leaf has been found, the proof becomes a not too difficult application of the maximum principle and of the compactification lemma in \cite{MV2} (see the end of section \ref{first}). Afterwards, in section \ref{second} we show how theorem \ref{main} applies to the case of some well-known foliations.

\

\noindent \textbf{Acknowledgements}: This paper was redacted while the author was a post-doc at the IMB in Dijon. I wish to thank L. Meersseman very much for his help, including originating the question and otherwise contributing to the article in many crucial ways.

\section{Main result}\label{first}
\subsection*{Statement}
Let $S$ be a $C^\infty$ Levi flat $2n+1$-manifold with boundary $C= bS$, and denote by $\mathcal F$ its smooth foliation by complex leaves. We will assume that $C$ is a compact complex manifold of dimension $n$.

We need to define a notion which extends that of partial compactification employed in \cite{MV2}.

\bd\label{exten}\emph{Let $M$ be a complex manifold of dimension $n$. We say that $M$  has an} end \emph{$E$ if there exists a sequence $\mathcal U_1\supset \mathcal U_2\supset\ldots$  such that every $b \oli{\mathcal U_j}$ is compact and connected and $\bigcap_i \oli{\mathcal U_i}=\emptyset$. Let now $X$ be another $n$-dimensional complex manifold, and let $\O\subset X$ be a proper subdomain such that $b\O$ is compact. We say that $X$} extends \emph{$M$ through its end $E$ if there exists a biholomorphism $\Phi:M\to \O$ such that $\bigcap_i \oli{\Phi(\mathcal U_i)}=b\O$.}
\emph{If $\O$ is dense in $X$ and $b\O\cong H$ is a compact $k$-dimensional complex submanifold of $X$, with $k<n$, we say - in accordance with the definition given in \cite{MV2} - that $X$ is a (partial)} holomorphic compactification \emph{of $M$ by $H$ at $E$-infinity.}

\emph{Assume, now, that $M=L$ is a leaf of a foliation $\mathcal F$ as before. Let $E=\{\mathcal U_i\}$ be an end of $L$, and suppose that (with respect to the topology of $S$) $\bigcap_i \oli{\mathcal U_i}=C$. In this situation, we say that $L$} ends at $C$ \emph{and an extension of $L$ at $E$-infinity is also said to be} at $C$-infinity.
\ed
Then, let $S$ be as above; we regard $C$ as a boundary leaf for $\mathcal F$. We are interested in a situation where
\begin{itemize}
\item [(A)] the (one sided) holonomy of $C$ is contracting and \textit{smoothly flat}, which means the following:  the holonomy group along $C$ takes values in the set of germs of real functions with infinite order of contact with the identity in $0$
\item [(B)] the holonomy covering $\widetilde C$ of $C$ extends through any of its ends.
\end{itemize}
Then we have
\bt \label{main} With the hypotheses above, there is no smooth embedding of $S$ as a Levi flat hypersurface (with boundary) of a complex manifold. \et
We remark that the hypotheses of theorem \ref{main} regard only the compact leaf $C$; in particular, we have that no embedding exists when $C$ satisfies some not too uncommon conditions, for example
\bc Let $S$ be a smooth Levi flat manifold whose foliation contains a compact complex leaf $C$ satisfying (A). If $C$ is a Riemann surface then $S$ admits no embeddings. \ec
We observe that if the holonomy of $C$ is tangent to identity but not smoothly flat, then there may exist an embedding of $S$. In fact, in \cite{Ba} Barrett shows an explicit construction of a (Lipschitz) embedding of $S^3$, with a foliation that is homeomorphic to Reeb's one but whose toric leaf's holonomy is not $C^\infty$ flat. In our context, since we are dealing with a situation with boundary, it is not difficult to give counterexamples where  the holonomy is even of class $C^k$ (see example \ref{Ck}).

\subsection*{Proof}

To prove theorem \ref{main}, as said before, we follow step by step the method employed by Barrett in \cite{Ba}. Assume, then, that there is a smooth embedding of $S$ into a complex $(n+1)$-manifold $M$; we will fix our attention to a neighborhood of the compact leaf $C$. We claim that
\bl \label{Ueda}There exists a holomorphic defining function $h$ for $C$, defined in a neighborhood of $C$ in $M$. Moreover, $h$ can be chosen in such a way that $d(Re h)|_S$ does not vanish in $C$. \el
To prove this lemma, we first give - following \cite{Ue} - a definition:
\bd \emph{Let $C$ be a compact complex hypersurface of a complex manifold $M$, and suppose that the normal bundle of $C$ is holomorphically trivial. Let $\mathfrak V=\{V_i\}$ be a small enough covering of a neighborhood of $C$ in $M$, and let $\mathfrak U = \{U_i\} = \{V_i\cap C\}$; then it is easy to see that there exists a system $\{w_i\}$ of local equations of $U_i$ in $V_i$ such that $w_i/w_k$ is well defined and equal to $1$ in $U_{ik}=U_i\cap U_k$. Denoting by $z_i$ a suitable set of local coordinates in $U_i$ (such that $(z_i,w_i)$ give coordinates for $V_i$), this means that for some integer $\nu$ and $f_{ik}\in \mathcal O(U_{ik})$ we have}
$$ w_k - w_i = f_{ik}(z_i)w_i^{\nu + 1} + o(\nu+1) $$
\emph{on $V_{ik}=V_i\cap V_k$. In such a case, the system $\{w_i\}$ is said to be} of type $\nu$. \emph{It is readily verified (see again \cite{Ue}) that $f_{ik}$ is a cocycle in $\mathcal Z^1(\mathfrak U,\mathcal O)$, and that it is a coboundary if and only if there exists a system of type $\nu + 1$. $C$ is said to be} of infinite type \emph{if any such system is a coboundary, i.e. there exists a system of type $\nu$ for all $\nu\in \mathbb N$.}
\ed
\br The type in the sense of Ueda defined above has the following geometrical meaning: it is the order of contact along $C$ of the line bundle $[C]$ (generated by $C$ as a divisor) and the trivial extension of the normal bundle to a neighborhood.\er
By hypothesis, the holonomy of the foliation of $S$ that we are considering along the compact leaf $C$ is trivial to infinite order. As a consequence of this fact, in the Appendix of \cite{Ba} the following is proven:

\bl The normal bundle of $C$ is holomorphically trivial; moreover, $C$ is of infinite type in $M$. \el

The point of the proof of the previous statement lies in the isomorphism (up to any finite order) between  the sheaf of functions which are locally constant on the leaves of $\mathcal F$ and a particular subsheaf of holomorphic functions of $M$. This isomorphism in turn depends on a result in \cite{BF} about the local (finite order) approximation of Levi flat hypersurfaces by zero sets of pluriharmonic functions, which holds for any dimension.

Lemma \ref{Ueda} is then a consequence of Theorem 3 in \cite{Ue}. Although that theorem is stated only for complex curves - since that is the framework of Ueda's paper - its proof works as well for any compact complex hypersurface of a complex manifold. In fact, the proof involves the construction of a new set of coordinate functions $\{u_i\}$ in $V_i$ which satisfy $u_i=u_j$ on $V_{ij}$. This is first carried out formally, expressing $u_i$ as a power series in $\{w_i\}$ with coefficients in $\mathcal O(U_i)$ in such a way that the relation is satisfied; the construction is possible because of the existence of a system of type $\nu$ for all $\nu\in \N$, which (roughly) implies the vanishing of the obstruction to the existence of each successive term of the series. The variables $z_i$ appear only through coefficients of the series in $\mathcal O(U_i)$, and the number of coordinates $z_i$ plays no role. Afterwards, the convergence of the series is estimated by a result by Kodaira and Spencer \cite{KS}, which is valid regardless of the dimension of $C$; although the argument is quite involved, the variables $z_i$ appear again only as a parameter.

\

\emph{Proof of \ref{main}}: Let $h$ be the function obtained by lemma \ref{Ueda}; $Re h$ has constant sign in a neighborhood of $C$ in $S$, we may suppose $Re h>0$. For a small enough $\e$, $\{0< Re h < \e\}$ is a one-sided tubular neighborhood $W$ of $C$ in $S$. A contradiction will be obtained by considering the behavior of the restriction of $h$ to $L\cap W$, where $L$  is a leaf in $S\setminus C$ whose closure contain $C$.
To this purpose, we first define a notion introduced in \cite{MV}, \cite{MV2}:

\bd \emph{we say that $\mathcal F$ is} tame \emph{if the following occurs: define the manifold $S'$  as
$$S' = S \sqcup (C\times [0,1]) / b S \sim C\times \{0\}$$
(i.e. $S'$ extends $S$ by attaching a collar $C \times (0,1]$ along $C$), and consider the foliation of $S'$ which agrees with $\mathcal F$ on $S$ and with the trivial one (induced by the submersion $C\times [0,1]\to [0,1]$) on $C\times [0,1]$. Moreover, endow the leaves of the foliation of $S'$ contained in $S$ with the complex structure inherited by $\mathcal F$, and each leaf contained in $C\times [0,1]$ with the complex structure of $C$. Then the foliation obtained is a smooth Levi foliation of $S'$. }\ed

For tame foliations, we can employ a compactification lemma proved in \cite{MV2}. Consider, on the restriction of $\mathcal F$ to a tubular neighborhood of $C$ in $S$, a second (tame) complex structure $J_1$ such that the structure induced on $C$ is the same as the original one. Then we can give the following variant of the compactification lemma cited above:

\bl\label{cmpctf} Let $L$ be a leaf of $\mathcal F$ ending in $C$, and $V$ a small enough tubular neighborhood of $C$. Let $L_1$ be the same leaf, but endowed by a complex structure $J_1$ as above. If $L_1$ admits an extension at $C$-infinity by a complex manifold $X$, then so does $L$. \el
\demo The proof employed in \cite{MV2} carries over: in fact, by the same argument the tameness of $\mathcal F$ implies that $J_1$ extends to $b\O$ smoothly (as an endomorphism of $T(X)$) and $J_1|_{b\O}= J|_{b\O}$. Hence $J_1$ extends smoothly on all of $X$, and it must be integrable since it is in $\O$ and $X\setminus \O$, so it is a complex structure in $X$.\enddemo

\noindent In order to apply lemma \ref{cmpctf} we first show - following the reparametrization method of Corollary 3 in  \cite{MV2} - that we can reduce to a tame situation.
Let $\widetilde C$ be the holonomy covering of $C$; then there exists an open subset $\widetilde V \subset \widetilde C \times [0,1)$, $C\times\{0\}\subset \widetilde V$, and a covering map $\pi$ from $\widetilde V$ onto a small tubular neighborhood $V$ of $C$ in $S$ which coincides with the holonomy covering on $\widetilde C\times \{0\}$ and sends each $(\widetilde C\times \{t\})|_{\widetilde V}$ diffeomorphically to a leaf $L_t$ of $\mathcal F|_V$. Let $\mathcal T_1,\ldots,\mathcal T_k$ be generators of the deck transformation group for $\pi$; then each $\mathcal T_j$ can be expressed as
 $$\mathcal T_j(p,t)= (T_j(p,t),d_j(t))$$
 where $p\in \widetilde C$ and $t\in [0,1)$. Now, consider the smooth automorphism of $\widetilde V$ defined by
 $$\Phi(p,t) = (p, \theta^{-1}(t)),\ \theta(t)=e^{\frac{1}{t}};$$
 then, if we endow $\widetilde V$ with the pull-back $CR$ structure, we obtain a manifold $\widetilde V'$ which is tautologically $CR$ isomorphic to $\widetilde V$. The quotient $V'$ of $\widetilde V'$ under the action of the group generated by
 $$\Phi^{-1}\circ \mathcal T_j\circ \Phi(p,t)= (T_j(p,\theta(t)),\theta^{-1}\circ d_j\circ \theta(t) )$$
is smoothly $CR$ isomorphic to $V$ and carries a tame foliation.

 Now, by (B) and lemma \ref{cmpctf} we deduce that a leaf $L$ ending in $C$ can be extended at $C$-infinity. In fact, if we endow each leaf $L$ of $V'$ with the structure $L^{pb}$ (obtained pulling back the complex structure of $C$) and we give to $\widetilde V'$ the trivial structure, the previously described covering $\widetilde V'\to V'$ is a biholomorphism along the leaves. It follows that each $L_t^{pb}$, hence $L_t$, can be extended. 

Consider, then, the biholomorphism $\Phi:L\to \O$ given by definition \ref{exten}; we are interested in $g=h\circ \Phi\in \mathcal O(\Phi^{-1}(L_W))$, where $L_W=L\cap W$. Since $h|_{L_W}$ converges to zero at the ending corresponding to $C$, we have that $g$ extends continuosly (by $0$) to $b\O$ and thus to $(X\setminus \O)\cup \Phi^{-1}(L_W)$. By Rado's theorem, then, follows that $g$ is holomorphic everywhere, hence $\O$ is actually dense  in $X$ and $b\O=H$ is an analytic subset of $X$. Then $Re h\circ \Phi$ is a non-constant pluriharmonic function on $\Phi^{-1}(L_W)\cup H$ which assumes minimum in its interior part (on $H$), a contradiction.

\section{Examples}\label{second}
\subsection*{Suspension of a Hopf manifold}

Fix coordinates $(z,w)$ in $\C^2$. As classified by Kodaira \cite{Ko}, any \emph{Hopf surface} is a quotient of $\C^2\setminus \{(0,0)\}$ by the action of
$$H:(z,w)\to (\a z+\l w^m, \b w)$$
where $m\in \N$ and $\a,\b,\l \in \C$ satisfy $(\b^m - \a)\l = 0$ and $0<|\a|\leq |\b| <1$.

Let, now, $\r:\R\to \R$ be a strictly increasing smooth function such that $\r(t)<t$ for $t>0$ and $\r(t) - t = o(t^d)$ as $t\to 0$ for all $d\in \N$. Consider $r: (\C^2\setminus \{(0,0)\}) \times \R \to (\C^2\setminus \{(0,0)\}) \times \R$ defined as
$$ r:(z,w,t) \to (H(z,w), \r(t)) $$
and let $S$ be the quotient of $(\C^2\setminus \{(0,0)\}) \times \R$ by the action of $r$. Since the action of $r$ preserves the foliation of $(\C^2\setminus \{(0,0)\}) \times \R$ by $\{t=const.\}$, $S$ inherits a structure of Levi flat manifold; the leaves are all isomorphic to $\C^2\setminus \{(0,0)\}$, except a compact leaf $C$ corresponding to $\{t=0\}$ which is a Hopf surface diffeomorphic to $S^3\times S^1$. Clearly, since $C$ is compact and non-Kähler, we know a priori that there is no embedding of $S$ as a Levi flat submanifold of either a Stein or a Kähler manifold. By \ref{main} we have that actually
\bc \label{Hopf} $S$ does not admit a $C^\infty$ embedding as a Levi flat hypersurface of a complex $3$-manifold.\ec
In this case, the holonomy covering of the Hopf surface $C$ coincides with its universal covering $\C^2\setminus \{(0,0)\}$, which has a partial holomorphic compactification by the $\C\P^1$ at infinity (and in fact the non-compact leaves are in turn compactifiable). Moreover, by the choice of $\r$ the holonomy of $C$ is contracting and $C^\infty$ flat, so that theorem \ref{main} applies.

On the other hand, in the non-smooth case the embedding is possible:

\bex\label{Ck} In fact, one can obtain an embedding in such a way that the holonomy of $C$ is flat up to any fixed order $d$: let $\r(t) = t - t^d$, and define the suspension as above. The resulting $S$ has a real analytic Levi foliation, and as such it can be embedded (see \cite{AF}). Notice that, since $\C^2\setminus\{(0,0)\}$ is compactifiable at both ends, the example also works for $\r(t)=t+t^d$.\eex

\subsection*{Partial generalization}
Let $P$ be a homogeneous polynomial in $\C^n$, and assume that $V=\{P=0\}$ is a smooth complex manifold outside the origin, with a smooth closure in $\C\P^n$. Choosing $0<\a<1$ and $\r:\R\to\R$ as above, we define the suspension
$$S = (V\setminus \{0\}) \times \R / \{(z,t) \sim (\a z, \r(t))\}.$$
We shall denote by $C$ the compact leaf, corresponding o $\{t=0\}$, of the foliation of $S$ induced by that of $(V\setminus \{0\}) \times \R$; the other leaves are isomorphic to $V\setminus \{0\}$.

As before, we have
\bc \label{gen} $S$ does not admit a $C^\infty$ embedding as a Levi flat hypersurface of a complex $n$-manifold.\ec
In this case, too, the holonomy covering of $C$ coincides with its universal covering; the partial compactification of $\widetilde C = V\setminus\{0\}$ is obtained by adding $\oli V \cap \C\P^{n-1}$, where $\oli V$ is the closure in $\C\P^n$.

\subsection*{Foliation of $S^5$}
In \cite{MV} it is constructed a smooth, $1$-codimensional Levi foliation of $S^5$, with two compact leaves. Each one of the compact leaves is isomorphic to a principal bundle over an elliptic curve $\mathbb E_\o$ whose fibers are elliptic curves. Since these compact leaves are not Kähler, it is once again clear that this foliation does not admit an embedding as Levi flat submanifold of a Stein or Kähler manifold. In fact
\bc \label{MV} There is no smooth embedding of $S^5$ as a Levi flat hypersurface of a complex manifold whose Levi foliation is diffeomorphic to the one obtained \cite{MV}.\ec
The foliation we are considering is constructed by gluing two partial ones. The first one lies in a tubular neighborhood $\mathcal N$ of $S^5\cap W$, where
$$W=\{(z_1,z_2,z_3)\in \C^3\setminus\{(0,0,0)\}: z_1^3+ z_2^3 + z_3^3=0 \}.$$
The foliation in $\mathcal N$ is obtained by considering a suitable covering
$$\C^\star\times B \to K\times D^2 \cong \mathcal N,$$ where $B\cong D^2\times S^1$ is a solid torus equipped with a foliation homeomorphic to Reeb's one. The compact leaf $S_\l$ lying in $\mathcal N$ is the restriction of the quotient to the boundary of $\C^\star\times B$. Thus, in a neighborhood of $S_\l$ the foliation is homeomorphic to the product of a disc by a neighborhood of the toric leaf in Reeb's foliation. In particular we have that the holonomy along $S_\l$ is trivial to infinite order.

In this case the holonomy covering of $S_\l$ does not coincide with its universal covering, but it is isomorphic to $W$; hence it admits a partial holomorphic compactification by $\oli W\cap \C\P^2_\infty$, where once again $\oli W$ is the closure in $\C\P^3$. It follows that theorem \ref{main} applies to this situation.

However, a more direct proof of the corollary can be achieved in the following way: let $h$ be as in lemma \ref{Ueda}. For a leaf $L$ sufficiently close to $S_\l$ and a small enough $\e$, the intersection
$$L\cap \{0<Re h<\e\}$$  is holomorphically equivalent to $D\times D^\star$, where $D$ is the unit disc and $D^\star$ is an annulus. The restriction of $Re h$ to $0\times D^\star$ is a positive harmonic function which vanishes, along with its conjugate, at $0$. But then $Re h$ extend to the whole disc, giving a contradiction by the maximum principle (see also \cite{Ba}).

\subsection*{Non-compactifiable cases}
We can give a simple example where $C$ admits no partial compactification, but hypothesis (B) of theorem \ref{main} applies, by the suspension of a compact complex curve of genus greater than $1$: consider such a curve $C$, let $\pi:D\to C$ be its universal covering and let $T_1,\ldots,T_k$ be the associated deck automorphisms. Let $S$ be the quotient of $D\times [0,1)$ by the action of the group generated by
$$ \mathcal T_1(z,t)=(T_1(z),r(t)),\ \mathcal T_j(z,t)=(T_j(z),t)\ {\rm for}\ j\geq 2 $$
where $r$ is chosen as in the previous examples. Then there is no embedding of $S$.  However, without respect to hypothesis (B), every time we can infer that the universal covering of the non-compact leaves is biholomorphic to a relatively compact domain of a complex manifold we can apply directly lemma \ref{Ueda} and conclude that no embedding exists.

One may conjecture that the flatness of the holonomy alone is sufficient to ensure that no embedding exists; the methods of section \ref{first}, though, are not sufficient to prove such a result.

\end{document}